\documentclass[12pt]{article}
\usepackage{amsthm,amsmath,amssymb}
\usepackage{a4wide}

\newtheorem{lemma}{Lemma}[section]

\begin{document}

\title{Two generalisations of the Binomial theorem}

\author{Sacha C. Blumen\footnotetext[1]{School of Mathematics and Statistics,
University of Sydney, NSW, 2006.  {\emph{E-mail:}} sachab@maths.usyd.edu.au}}
\date{ }

\maketitle

\begin{abstract}
We prove two generalisations of the Binomial theorem that are also
generalisations of the $q$-binomial theorem.  These generalisations arise from
the commutation relations satisfied by the components of the 
co-multiplications of non-simple root vectors in the quantum superalgebra 
$U_{q}(osp(1|2n))$.
\end{abstract}

Generalisations of the Binomial theorem can be used to expand
powers of sums of elements of non-abelian associative algebras.
One well-known generalisation is the $q$-binomial theorem, which gives the expansion of
$(x+y)^{n}$, for each $n = 1, 2, 3, \ldots$, 
where $x$ and $y$ are non-commuting quantities 
satisfying $xy = q yx$, for some $0 \neq q  \in \mathbb{C}$ 
also satisfying $q^{2} \neq 1$.

Quantum algebras and quantum superalgebras are 
a rich source of elements satisfying the relation $xy = q yx$.  
In addition, many elements of quantum algebras and quantum superalgebras 
satisfy much more complicated relations, 
and different generalisations of the Binomial theorem may be called on in performing
 calculations in these or other algebras.

The two generalisations of the Binomial theorem in this note  
appear in the author's Ph.D thesis \cite{blumen2005}.  
For readers familiar with quantum algebras and quantum superalgebras, 
these two generalisations  are related to the commutation relations
satisfied by the components of the co-multiplications of non-simple root vectors in $U_{q}(osp(1|2n))$
 defined following \cite{kt}. 
Using these generalisations, it was shown in \cite{blumen2005} that
a certain two-sided ideal $I$ 
of $U_{q}(osp(1|2n))$ is also a Hopf ideal 
when $q = \exp{(2 \pi i/N)}$ for some integer $N \geq 3$,  a
consequence of which is that the quotient algebra $U_{q}^{(N)}(osp(1|2n)) = U_{q}(osp(1|2n))/I$ 
admits a universal $R$-matrix originally written down in \cite{zhang1992}.

The results in this note may be
useful in calculations in other quantum (super)algebras, 
but I leave this for further exploration.

\begin{section}{Notations}

We write $\mathbb{N} = \{1, 2, 3, \ldots\}$ and 
$\mathbb{Z}_{+} = \mathbb{N} \cup \{0\}$.
For all $0 \neq q  \in \mathbb{C}$ satisfying $q^{2} \neq 1$, we define
$$[n]_{q} = \frac{1-q^{n}}{1-q}, \hspace{10mm}
  [n]_{q}! = [n]_{q} [n-1]_{q} \cdots [1]_{q}, \hspace{10mm} 
  [0]_{q}! = 1, \hspace{10mm}
  \mbox{for each } n \in \mathbb{Z}_{+}.$$

\end{section}

\begin{section}{Two generalisations of the Binomial theorem}

\begin{lemma}
\label{lem:generalisation_A}
Let $a, b$ and $c$ be elements of an associative algebra over $\mathbb{C}$ satisfying 
\begin{equation}
\label{eq:apple(1)}
ab = q ba + c, \hspace{10mm} ac = q^{2} ca, \hspace{10mm} cb = q^{2} bc,
\end{equation}
where $0 \neq q \in \mathbb{C}$ and $q^{2} \neq 1$, then
$$(a + b)^{n} = \sum_{\stackrel{\alpha,\beta,\gamma \in \mathbb{Z}_{+}}{\alpha + 2\beta + \gamma = n}}  
\frac{ [n]_{q}!  }{ [\alpha]_{q}! [\gamma]_{q}! [2]_{q} [4]_{q} \cdots [2\beta]_{q} }   \ 
b^{\alpha}c^{\beta} a^{\gamma}, \hspace{10mm} n \in \mathbb{N}.$$
\end{lemma}
\begin{proof}
From (\ref{eq:apple(1)}) we can inductively prove that
\begin{equation}
\label{eq:apple(2)}
a^{n}b = q^{n}  b a^{n} + q^{n-1}[n]_{q}  c a^{n-1}, \hspace{10mm} n \in \mathbb{N},
\end{equation}
and we can use (\ref{eq:apple(2)}) to prove the following relation, 
where we fix $\alpha, \beta, \gamma$ to be non-negative integers:
\begin{eqnarray*}
b^{\alpha}c^{\beta} a^{\gamma} b & = & 
q^{\gamma + 2\beta} \ b^{\alpha + 1} c^{\beta} a^{\gamma} + 
q^{\gamma-1} [\gamma]_{q} \ b^{\alpha} c^{\beta + 1} a^{\gamma-1}.
\end{eqnarray*}
It is not dificult to show that $\alpha + 2\beta + \gamma = n$
if $b^{\alpha}c^{\beta} a^{\gamma}$ is a component in the expansion of $(a+b)^{n}$, 
thus we have
\begin{equation}
\label{appendixB:simoneyoung}
(a + b)^{n} = \sum_{\stackrel{\alpha, \beta, \gamma \in \mathbb{Z}_{+}}{\alpha + 2\beta + \gamma = n}}  
\theta(\alpha,\beta,\gamma) \ b^{\alpha}c^{\beta} a^{\gamma}, \hspace{10mm} n \in \mathbb{N}, 
\end{equation}
for some set of coefficients 
$\left\{\theta(\alpha,\beta,\gamma) \in \mathbb{C} | \ \alpha, \beta, \gamma \in \mathbb{Z}_{+} \right\}$.

From the algebra relations and (\ref{appendixB:simoneyoung}), 
the coefficients $\theta(\alpha,\beta,\gamma)$ satisfy the recurrence relation:
\begin{eqnarray}
\theta(\alpha,\beta,\gamma)
& = & \theta(\alpha,\beta,\gamma-1) + q^{\gamma + 2\beta} \theta(\alpha-1,\beta,\gamma)
+ q^{\gamma}[\gamma+1]_{q} \theta(\alpha,\beta-1,\gamma+1)  \nonumber \\
& & \label{eq:appenB:binom1}
\end{eqnarray}
subject to the boundary conditions $\theta(1,0,0)=\theta(0,0,1)=1$.
In (\ref{eq:appenB:binom1}) we fix $\theta(\alpha,\beta,\gamma) = 0$ 
if any of $\alpha, \beta$ or $\gamma$ are negative.
To complete the proof, we just need to show that the set of constants 
\begin{equation}
\label{appendixB:simoneyoung2}
\theta(\alpha, \beta, \gamma) = 
\frac{[\alpha + 2\beta + \gamma]_{q}!}{[\alpha]_{q}![\gamma]_{q}![2]_{q}[4]_{q} \cdots [2\beta]_{q}},
\end{equation}
solves the recurrence relation and also satisfies the boundary conditions.  
The latter is easy to see, and substituting (\ref{appendixB:simoneyoung2}) 
into the right hand side of (\ref{eq:appenB:binom1}) gives 
\begin{eqnarray*}
& & \frac{[\alpha + 2\beta + \gamma - 1]_{q}!}
{[\alpha]_{q}! [\gamma]_{q}! [2]_{q} [4]_{q} \cdots [2\beta]_{q}} 
\left( [\gamma]_{q} + q^{\gamma + 2\beta}[\alpha]_{q} + q^{\gamma}[2 \beta]_{q}  \right)  
\end{eqnarray*}
which equals the right hand side of (\ref{appendixB:simoneyoung2}) as desired.
\end{proof}
Note that Lemma \ref{lem:generalisation_A} is just the $q$-binomial theorem
when $c=0$ in (\ref{eq:apple(1)}).
For readers familiar with quantum superalgebras,
the generalisation of the Binomial theorem in 
Lemma \ref{lem:generalisation_A} is connected to the relations satisfied by the components 
of  $\Delta(e_{\zeta}) \in U_{q}(osp(1|2n))$ where
$\zeta$ is a non-simple root containing {\emph{one}} copy of the odd simple root of
$osp(1|2n)$.
Similar remarks apply for Lemma \ref{lem:generalisation_B}, but here
$\zeta$ is a non-simple root containing {\emph{two copies}} of the odd simple root.

\begin{lemma}
\label{lem:generalisation_B}
Let $a, b$ and $c$ be elements of an associative algebra over $\mathbb{C}$ satisfying 
\begin{equation}
\label{lem2:relations}
ac = q^{2} ca + \xi b^{2}, \hspace{10mm} ab = q^{2} ba, \hspace{10mm} bc = q^{2} cb,
\end{equation}
where $0 \neq q \in \mathbb{C}$, $q^{2} \neq 1$ and $\xi =  -(1+q)^{2}/(q-q^{-1})$, then
$$(a + b + c)^{n} = \sum_{\stackrel{\alpha,\beta,\gamma \in \mathbb{Z}_{+}}{\alpha + \beta + \gamma = n}} 
\frac{ [n]_{q^{2}}! \ \phi_{\beta} }{ [\alpha]_{q^{2}}! [\beta]_{q^{2}}! [\gamma]_{q^{2}}! } \ 
c^{\alpha} b^{\beta} a^{\gamma}, \hspace{10mm}  n \in \mathbb{N},$$
where $\phi_{\beta} \in \mathbb{C}$ is recursively defined by
$$\phi_{0}=1, \hspace{10mm}  \phi_{1} = 1,  
\hspace{10mm} \phi_{\beta} = \phi_{\beta-1} + \xi [\beta-1]_{q^{2}} \phi_{\beta-2},
\hspace{10mm} \beta \in \mathbb{N} \backslash \{1\}.$$

\end{lemma}
\begin{proof}
From (\ref{lem2:relations}) we can inductively prove that
$$a^{n}c = q^{2n}ca^{n} + \xi q^{2(n-1)}[n]_{q^{2}} b^{2} a^{n-1}, 
\hspace{10mm} n \in \mathbb{N},$$
whch we can use to show the following relations, 
where $\alpha, \beta$ and $\gamma$ are all non-negative integers:
\begin{eqnarray*}
c^{\alpha} b^{\beta} a^{\gamma} b & = & q^{2\gamma} c^{\alpha} b^{\beta+1} a^{\gamma} \\
c^{\alpha} b^{\beta} a^{\gamma} c & = & q^{2\gamma + 2\beta} c^{\alpha+1} b^{\beta} a^{\gamma}
+ \xi q^{2(\gamma-1)}[\gamma]_{q^{2}} c^{\alpha} b^{\beta+2} a^{\gamma-1}.
\end{eqnarray*}
It is not difficult to show that 
$\alpha + \beta + \gamma = n$ if $c^{\alpha} b^{\beta} a^{\gamma}$ is a component of 
$(a+b+c)^{n}$,  thus we have
\begin{equation}
\label{appendixB:simoneyoung3}
(a+b+c)^{n} = \sum_{\stackrel{\alpha, \beta, \gamma \in \mathbb{Z}_{+}}{\alpha + \beta + \gamma = n}} 
\theta(\alpha,\beta,\gamma) \ c^{\alpha} b^{\beta} a^{\gamma}, \ 
\hspace{10mm} n \in \mathbb{N},
\end{equation}
for some set of coefficients 
$\{ \theta(\alpha,\beta,\gamma) \in \mathbb{C} | \ \alpha,\beta,\gamma \in \mathbb{Z}_{+} \}$.

From (\ref{appendixB:simoneyoung3}) and the algebra relations,
the coefficients $\theta(\alpha,\beta,\gamma)$ satisfy the following recursion relation
\begin{eqnarray}
\theta(\alpha,\beta,\gamma) & = & \theta(\alpha,\beta,\gamma-1) + q^{2\gamma} \theta(\alpha,\beta-1,\gamma)
+ q^{2\gamma + 2\beta} \theta(\alpha-1,\beta,\gamma)  \nonumber  \\ 
& & + \xi q^{2\gamma} [\gamma+1]_{q^{2}} \theta(\alpha,\beta-2,\gamma+1) \label{eq:appenB:binom2}
\end{eqnarray}
subject to the boundary conditions $\theta(1,0,0)= \theta(0,1,0) = \theta(0,0,1) = 1$.
In (\ref{eq:appenB:binom2}) we fix $\theta(\alpha,\beta,\gamma)=0$ if any of $\alpha, \beta, \gamma$ are negative.
To complete the proof, we just need to show that the set of constants
\begin{equation}
\label{appendixB:simoneyoung4}
\theta(\alpha,\beta,\gamma) = 
\frac{ [\alpha + \beta + \gamma]^{q^{2}}! \ \phi_{\beta} }
{ [\alpha]^{q^{2}}! [\beta]^{q^{2}}! [\gamma]^{q^{2}}! },
\end{equation}
solves the recurrence relation and satisfies the boundary conditions,
where $\phi_{\beta}$ is itself recursively defined as stated in the lemma.

It is clear that the constants in
(\ref{appendixB:simoneyoung4}) satisfy the boundary conditions, 
and substituting them into the right hand side of (\ref{eq:appenB:binom2}) gives 
\begin{eqnarray}
\lefteqn{  
\frac{ [\alpha+\beta+\gamma-1]_{q^{2}}! \ \phi_{\beta} }{ [\alpha]_{q^{2}}! [\beta]_{q^{2}}! [\gamma-1]_{q^{2}}! }
+ 
q^{2\gamma}\frac{ [\alpha+\beta+\gamma-1]_{q^{2}}! \
\phi_{\beta-1}}{[\alpha]_{q^{2}}![\beta-1]_{q^{2}}![\gamma]_{q^{2}}! }
 } \nonumber \\
 & & + q^{2\gamma + 2\beta} 
\frac{ [\alpha+\beta+\gamma-1]_{q^{2}}! \ \phi_{\beta} }{ [\alpha-1]_{q^{2}}! [\beta]_{q^{2}}! [\gamma]_{q^{2}}! }
 + \xi q^{2\gamma}[\gamma+1]_{q^{2}}
 \frac{ [\alpha+\beta+\gamma-1]_{q^{2}}! \ \phi_{\beta-2} }{ [\alpha]_{q^{2}}! [\beta-2]_{q^{2}}!
 [\gamma+1]_{q^{2}}! } \nonumber \\
 & = & \frac{ [\alpha+\beta+\gamma-1]_{q^{2}}!  }{ [\alpha]_{q^{2}}! [\beta]_{q^{2}}! [\gamma]_{q^{2}}! }
 \left( [\gamma]_{q^{2}} \phi_{\beta} + q^{2\gamma + 2\beta}[\alpha]_{q^{2}}\phi_{\beta}
 + q^{2\gamma} [\beta]_{q^{2}} \left[\phi_{\beta-1} + \xi[\beta-1]_{q^{2}}\phi_{\beta-2}\right]\right).
 \nonumber \\
 & &  \label{appendixB:simoneyoung5}
\end{eqnarray}
By writing $\phi_{\beta} = \phi_{\beta-1} + \xi[\beta-1]_{q^{2}} \phi_{\beta-2}$
for each $\beta \in \mathbb{N} \backslash \{1\}$,
 we can rewrite
(\ref{appendixB:simoneyoung5}) as
$$\frac{ [\alpha+\beta+\gamma-1]_{q^{2}}!  }{ [\alpha]_{q^{2}}! [\beta]_{q^{2}}! [\gamma]_{q^{2}}! }
 \left( [\gamma]_{q^{2}} \phi_{\beta} + q^{2\gamma + 2\beta}[\alpha]_{q^{2}}\phi_{\beta}
 + q^{2\gamma} [\beta]_{q^{2}} \phi_{\beta} \right),$$
which is just 
$\displaystyle{\frac{ [\alpha+\beta+\gamma]_{q^{2}}! \ \phi_{\beta}}{[\alpha]_{q^{2}}!
[\beta]_{q^{2}}! [\gamma]_{q^{2}}!}}$
as claimed.
\end{proof}

Note that Lemma \ref{lem:generalisation_B} is just a version of the $q$-multinomial theorem if we
artificially fix $\xi=0$ in  (\ref{lem2:relations}).
We obtain a general expression for $\phi_{\beta}$ below.
\begin{lemma}
\label{lem:calculation_of_phi}
Let $0 \neq q \in \mathbb{C}$ satisfy $q^{2} \neq 1$ and let 
$\phi_{\beta} \in \mathbb{C}$ be recursively defined by
$$\phi_{0}=1, \hspace{10mm}  \phi_{1} = 1,  \hspace{10mm} 
\phi_{\beta} = \phi_{\beta-1} + \xi [\beta-1]_{q^{2}} \phi_{\beta-2},
\hspace{10mm} \beta \in \mathbb{N} \backslash \{1\},$$
where $\xi = -(1+q)^{2}/(q-q^{-1})$.  
Then $\phi_{\beta}$ is explicitly
$$\phi_{0} = 1, \hspace{10mm} \phi_{1}=1, \hspace{10mm}
\phi_{2i}   = (1-q)^{-i} \Psi_{2i},  \hspace{10mm} \phi_{2i+1} = [2i+1]_{q} \phi_{2i},$$
for each $i \in \mathbb{N}$, where
$$\Psi_{2i} = \frac{[4]_{q}}{[2]_{q}} [3]_{q} \frac{[8]_{q}}{[4]_{q}} [5]_{q} 
\frac{[12]_{q}}{[6]_{q}} [7]_{q} \cdots [2i-1]_{q} \frac{[4i]_{q}}{[2i]_{q}}.$$
\end{lemma}
\begin{proof}
We firstly calculate that
$$
\phi_{2}= 1 + \xi [1]_{q^{2}}  = (1 + q^{2})/(1-q) = (1-q)^{-1} [4]_{q}/[2]_{q},
$$
%
and thus $\phi_{\beta}$ is as claimed for $\beta=0,1,2$.
Assume now that $\phi_{2i}$ and $\phi_{2i-1}$ are as is stated in the lemma for some $i \in \mathbb{N}$, 
then 
\begin{eqnarray}
\phi_{2i+1} & = & \phi_{2i} - \frac{(1+q)^{2}}{q-q^{-1}}[2i]_{q^{2}}\phi_{2i-1} \nonumber \\
            & = & \frac{(1+q)\left(-q^{-1} - [2i]_{q} \right)}{q-q^{-1}} [2i-1]_{q}
	    \frac{[4i]_{q}}{[2i]_{q}} \phi_{2i-2} \nonumber \\
	    & = & (1-q)^{-1} [2i-1]_{q} \frac{[4i]_{q}}{[2i]_{q}} [2i+1]_{q} \phi_{2i-2} \nonumber \\
	    & = & [2i+1]_{q} \phi_{2i}, \hspace{5mm} \mbox{and} \nonumber \\
\phi_{2i+2} & = & [2i+1]_{q}\phi_{2i}  - \frac{(1+q)^{2}}{q-q^{-1}}[2i+1]_{q^{2}} \phi_{2i} \nonumber \\
            & = & \left([2i+1]_{q} - \frac{(1+q)}{q-q^{-1}}[4i+2]_{q}  \right)\phi_{2i} 
	           \label{appendixB:simoneyoung7} \\
	    & = &  \frac{-q^{-1}(1+q) [2i+1]_{q}[4i+4]_{q}}{ (q-q^{-1}) [2i+2]_{q} } \phi_{2i}  
	           \nonumber \\
	    & = & (1-q)^{-1} [2i+1]_{q} \frac{[4i+4]_{q}}{[2i+2]_{q}} \phi_{2i},  \nonumber
\end{eqnarray}
as claimed.
We used $(1+q)\left[2i+1\right]_{q^{2}} = \left[4i+2\right]_{q}$ 
to obtain (\ref{appendixB:simoneyoung7}).

\end{proof}

The reader may find it interesting to explore these two generalisations of the Binomial theorem when
$q=\pm\exp{(2 \pi i/N)}$ and $N \geq 3$ an integer.

\end{section}

\noindent
{\bf{Acknowledgments}}

\noindent
I would like to thank Dr Chris Cosgrove, who helped in a related problem
and whose thoughts subsequently led to a useful idea in this work.

\end{document}